\documentclass{amsart}

\usepackage{amssymb}
\usepackage{amscd}
\usepackage{amsmath}
\usepackage{amsfonts}
\usepackage{a4}
\usepackage{amsthm}
\usepackage{color}

\theoremstyle{plain}

\newtheorem{theo}{Theorem}[section] 
\newtheorem{prop}[theo]{Proposition}
\newtheorem{lemme}[theo]{Lemme}

\theoremstyle{definition}

\newcommand{\eval}[1]{\lvert_{#1}}

\newcommand{\Ci}{C^{\infty}}

\newcommand{\R}{{\mathbb{R}}}

\newcommand{\C}{{\mathbb{C}}}
\newcommand{\Z}{{\mathbb{Z}}}

\newcommand{\de}{{\delta}}
\newcommand{\be}{{\beta}}
\newcommand{\al}{{\alpha}}
\newcommand{\la}{{\lambda}}
\newcommand{\De}{{\Delta}}
\newcommand{\si}{{\sigma}}
\newcommand{\hb}{{\hbar}}

\newcommand{\ga}{{\gamma}}
\newcommand{\te}{{\theta}}
\newcommand{\om}{{\omega}}
\newcommand{\Om}{\Omega}
\newcommand{\ph}{{\varphi}}
\newcommand{\Ph}{{\Phi}}
\newcommand{\Ga}{{\Gamma}}

\newcommand{\La}{{\Lambda}}
\newcommand{\Si}{{\Sigma}}

\newcommand{\operateur}{{\mathcal{Q}}}
\newcommand{\identite}{{\operatorname{Id}}}
\newcommand{\Trace}{{\operatorname{Tr}}}
\newcommand{\trace}{{\operatorname{tr}}}
\newcommand{\volume}{{\operatorname{Vol}}}

\newcommand{\Hilbert}{{\mathcal{H}}}

\newcommand{\Toeplitz}{{\mathcal{T}}}
\newcommand{\Lie}{{\mathcal{L}}}
\newcommand{\equivalence}{{\mathcal{E}}}

\newcommand{\contravariant}{{\operatorname{cont}}}  
\newcommand{\normalised}{{\operatorname{norm}}}

\newcommand{\Lag}{{\operatorname{Lag}}}

\newcommand{\Mas}{{\mathcal{M}}}

\newcommand{\ideal}{{\mathcal{I}}}
\newcommand{\Sym}{\operatorname{Sym}}
\newcommand{\no}{{\mathcal{N}}}

\newcommand{\scratch}[1]{}

\author{L. Charles} 

\address{Universit{\'e} Pierre et Marie Curie-Paris6, UMR 7586 Institut de
  Math{\'e}matiques de Jussieu, Paris, F-75005 France.}

\email{charles@math.jussieu.fr}



\title{Symbolic calculus for Toeplitz operators with half-forms} 

\begin{document}

\begin{abstract}
This paper is devoted to the use of half-form bundles in the symbolic
calculus of Berezin-Toeplitz operators on K{\"a}hler manifolds. We state
the Bohr-Sommerfeld conditions and relate them to the
functional calculus of Toeplitz operators, a trace formula and the
characteristic classes in deformation quantization. We also develop
the symbolic calculus of Lagrangian sections, with the crucial
estimates of the subprincipal terms. 
\end{abstract}

\maketitle

\bibliographystyle{plain}

In semi-classical analysis we usually consider (pseudo) differential
operators depending on a small parameter and acting on a $L^2$ space,
the underlying classical limit being a cotangent space with its
canonical symplectic structure. In this paper we are interested in a
similar theory where the classical phase space is a compact K{\"a}hler
manifold endowed with a prequantum line bundle $L$. Here the quantum
Hilbert space consists of the holomorphic sections of $L^k$. The
small parameter is the inverse of the power $k$. The operators of
interest are the
Berezin-Toeplitz operators. This setting was mainly introduced by
Kostant \cite{Ko},
Souriau \cite{So} and Berezin \cite{Be} and the suitable microlocal techniques were
developed by Boutet de Monvel and Guillemin \cite{BoGu}. Since then many standard
results for pseudo-differential operators have been adapted to this
context, like for instance the Schnirelman theorem \cite{Ze}, the Gutzwiller trace formula
\cite{BP2} or 
the Bohr-Sommerfeld conditions \cite{oim2}. The statements of these results are
easily predictable as far as only the symplectic structure of the
phase space is concerned, because they are the same
for the cotangent and K{\"a}hler spaces. But the semi-classical results
for the pseudo-differential operators may involve also some invariants,
like the subprincipal symbol or the Maslov index, which
do not only depend on the symplectic structure and consequently 
are difficult to identify in the K{\"a}hler setting. Furthermore these
quantities generally appear as quantum corrections and are
difficult to compute. Nevertheless in the papers \cite{oim1} and \cite{oim2}, we
carried out successfully some techniques to handle this. To formulate our
result we used the Riemannian metric of the K{\"a}hler structure instead
of the vertical polarization of the cotangent bundle. Typically we
proved some Bohr-Sommerfeld conditions where the Maslov index is
replaced with a curvature integral. Actually we missed the right
formulation which uses the half-form bundles. 
The main purpose of this paper and the sequel \cite{new} is to develop this point
of view. In this part we focus on the Bohr-Sommerfeld conditions
whereas \cite{new} is devoted to the dependence of the quantization on the complex structure.

Concretely we alter the usual setting by defining the quantum space  as the
space of holomorphic sections of $L^k \otimes L_1 \otimes \delta
\rightarrow M$. Here
$L$ is the prequantum bundle as previously, $L_1$ is an auxiliary
line bundle and $\delta$ is a half-form bundle, {\em i.e.} a square
root of the canonical bundle of $M$. A priori artificial, this
decomposition enlightens the semi-classical results, even in the usual
case where $L_1 \otimes
\delta$ is the trivial bundle. Roughly speaking the contribution of
$L_1$ in the semi-classical limit is the same as the one
of $L^k$ (one can view $L^k \otimes L_1$ as a first order
deformation of $L^k$) whereas the half-form bundle contributes in a
specific way. This principle will be confirmed in all our
results.  Another important point is that there
is a topological obstruction to the existence of half-form bundles. To
avoid this problem we consider globally the bundle $L^k \otimes K$ and
write locally $K = L_1 \otimes \delta$. 
The situation here is analogous to that in Riemannian geometry where we
think any Clifford module, at least locally, as the spinor bundle
twisted with an auxiliary bundle.

The first section is devoted to basic properties of Toeplitz
operators and their symbolic calculus. In particular an important
subprincipal symbol is defined. We state the Bohr-Sommerfeld
conditions in section 2 and relate them to the symbolic calculus and
trace formula by adapting an argument of Colin de Verdi{\`e}re
\cite{Co}. Here the formulation with half-forms permits to check easily the
consistency of the results. The next sections contain the proof of the
Bohr-Sommerfeld conditions. In section 4, we introduce the Lagrangian
sections, which are similar to the Lagrangian distributions, and
develop their symbolic calculus. Bohr-Sommerfeld conditions follows immediately. 
A comparison with the usual setting
is included, where a $\Z_4$-bundle plays a role analogous to the
Maslov bundle. In section 5, the technical part of the paper, we
provide the proof for the symbolic calculus of  the Lagrangian
sections. We follow essentially the method of \cite{oim2} but avoid the
complicated computations involving the derivatives of the K{\"a}hler
metric. These simplifications rely on a version of the stationary
phase lemma stated in an appendix of this paper. In view of this
proof, we think that our result  should generalize mutatis mutandis to the case where the
symplectic manifold doesn't admit any integrable complex structure.

{\bf Acknowledgments} We thank Y. Colin de Verdi{\`e}re who provided us
his preprint \cite{Co} and suggested us to adapt his argument to the Toeplitz
operators. This was actually one of our original motivations to develop
the half-form formalism. We also thank F. Faure for his kind interest.

\section{The setting}
\subsection{Square root of line bundle}

Let $M$ be a manifold and $F \rightarrow M$ be a complex line bundle.  A {\em square root} $(\delta,
\ph)$ of $F$ is a line bundle $\delta \rightarrow M$ together with an
isomorphism of line bundle $\ph : \delta^{\otimes 2} \rightarrow
F$. If $M$ is a complex manifold, a square root of its canonical bundle
$\La ^{n,0} T^*M$ is called a {\em half-form bundle}. 
Let us state basic properties of square roots. 

If $F$ has a Hermitian structure and $(\delta,
\ph)$ is a square root of $F$, then $\delta$ has a unique Hermitian structure
such that $\ph$ is a isomorphism of Hermitian line bundle. In the same
way, if $F$ is holomorphic or flat, $\delta$ inherits the same
structure. If $D^F$ is a first order differential operator acting on
sections of $F \rightarrow M$, then there exists a unique first
order differential operator $D^\delta$ acting on section of $\delta$ such that 
$$ D^F \ph (s \otimes s ) = 2 \ph ( s \otimes D^\delta s), \quad \forall \; s \in \Ci
(M,\delta).$$
A line bundle admits a square root if and only if its Chern class is divisible by 2 in $H^2(M, \Z)$. 
Two square
roots $(\delta, \ph)$ and $(\delta' , \ph')$ of $F$ are {\em equivalent}
if there exists an isomorphism $\Psi : \delta \rightarrow
\delta'$ such that $\ph' \circ \Psi^2 = \ph$. 

\begin{prop} Assume that $F$ admits a square root. Then the set
  of
  equivalence classes of square roots of $F$ is a principal
  homogeneous space for the first group of cohomology of $M$ with coefficient
  in $\Z_2$. 
\end{prop}

\begin{proof} 
First, if $(\delta, \ph)$ is a square root of the trivial line bundle
$1_M = M \times \C$, then $\delta$ inherits a flat structure from $1_M$
with structure group $\Z_2$. Furthermore this flat structure
determines $\ph$. It is easily proved that this induces an
isomorphism between the set of equivalence classes of square roots of
$1_M$ and the set of equivalences of flat line bundles with structure
group $\Z_2$. The latter is isomorphic to $H^1(M, \Z_2)$. Now, observe
that the tensor product of a square root of $L$ with a square root of
$1_M$ is a square root of $L$. This defines an action of $H^1(M,
\Z_2)$ on  $\equivalence _L$, which is easily shown to be free and transitive.
\end{proof}

\subsection{Quantum spaces} 
Let $M$ be a connected compact K{\"a}hler manifold of complex dimension
$n$. Denote by $\om \in \Om ^2 (M, \R)$ the fundamental form of $M$. Assume $M$ is endowed with a prequantization bundle $$L
\rightarrow M,$$ that is a
Hermitian line bundle with a connection $\nabla^L$ of curvature $\frac{1}{i}
\om$. Since $\om$ is a $(1,1)$-form, $L$ has a natural holomorphic
structure defined in such a way that the (local) holomorphic sections
satisfy the Cauchy-Riemann equations: $\nabla _{\bar{Z}} s
= 0$ for every holomorphic vector field $Z$ of $M$.

Let $K \rightarrow M$ be a Hermitian holomorphic line bundle. For
every positive integer $k$ define the quantum space
$\Hilbert_k$ : 
$$ \Hilbert_k = \bigl\{ \text{holomorphic section of } L^k \otimes K
\bigr\}.$$  
Assume that $M$ carries a half-form bundle $(\delta, \varphi)$.
$\delta \rightarrow M$ inherits a Hermitian scalar product and a holomorphic structure from
$\La^{ n, 0 }T^* M$. Introduce the Hermitian holomorphic line bundle $L_1$
such that $$K = L_1 \otimes \delta$$
and let $\frac{1}{i} \om_1 $ be the curvature of the Chern connection of
$L_1$.

Since $M$ is
compact, $\Hilbert_k$ is finite dimensional and it follows
from  the Riemann-Roch-Hirzebruch theorem and
Kodaira vanishing theorem that  
\begin{gather} \label{Riemann-Roch} 
 \operatorname{dim} \Hilbert_k = \Bigl( \frac{k}{2 \pi} \Bigr)^{n}
\int_M ( \om + k^{-1} \om_1 )^{\wedge n} /n!  + O(k^{n-2}) \end{gather} 
To interpret this formula, we consider $L^k \otimes L_1$ and $\om +
\hbar \om_1$ as
deformations of $L^k$ and $\om$ which give the first quantum
corrections in the semi-classical limit. Indeed  the leading term 
$$\bigl( \tfrac{k}{2 \pi} \bigr)^{n} \int \om^{\wedge n} /n!$$  gives
the second-order correction when we replace $\om$ with $\om + k^{-1}
\om_1$. Furthermore in the case $M$ doesn't carry any half-form
bundle, equation \eqref{Riemann-Roch}  is still
valid if we define $\om_1$ by 
$$ \om_1 :=  \om_K - \om_c/2,$$ 
where $\frac{1}{i} \om_K$
and $\frac{1}{i} \om_c $ are the curvatures of the Chern connections of
$K$ and $\La^{n,0} T^*M$.

\subsection{Toeplitz operators} 
Denote by $\Pi_k$ the orthogonal projector of $L^2(M,
L^k \otimes K)$ onto ${\mathcal{H}}_k$, where the scalar product of two sections of $L^k \otimes K$ is
defined from the Hermitian structures of $L$ and $K$ and the
Liouville form $\mu_M$.

A Toeplitz operator is any  sequence $(T_k: \Hilbert_k \rightarrow
\Hilbert_k )$ of operators of the form 
$$ T_k = \Pi_k f(.,k) + R_k , $$
where $f(.,k)$ is a sequence of $\Ci( M)$ with an asymptotic expansion
$f_0 + k^{-1} f_1 +...$ for the $\Ci$ topology and the norm of $R_k$
is $O(k^{-\infty})$. 

The set $\Toeplitz$  of Toeplitz operators is a semi-classical algebra associated
to $(M, \om)$ in the following sense. 

\begin{theo} \label{sc_algebra}
$\Toeplitz$ is closed under the
formation of product. So it is a star algebra, the identity is
$(\Pi_k)$. The symbol
map $$\si_{ \contravariant} : \Toeplitz \rightarrow \Ci(M)[[\hb]],$$ sending $T_k$ into the
formal series $f_0 + \hb f_1 + ...$ where the functions $f_i$ are the
coefficients of  the asymptotic expansion of the multiplicator
$f(.,k)$, is
well defined.  It is onto and its kernel is the ideal consisting of
$O(k^{-\infty})$ 
Toeplitz operators. More precisely for any integer $\ell$, 
$$ \| T_k \| = O(k^{-\ell} ) \text{ if and only if } \si_{\contravariant} (T_k)
= O( \hb^{\ell}). $$ 
 Furthermore, the induced product $*_{\contravariant}$ on $\Ci (M) [[\hbar]]$ is a
 star-product. 
\end{theo}
Following the terminology of Berezin in \cite{Be}, we call $\si_{
  \contravariant}$  the
contravariant symbol map. This result is essentially a consequence of the works of Boutet de Monvel and Guillemin
\cite{BoGu} (cf. also \cite{Gu}, \cite{BoMeSc} and \cite{oim1}). 
Let us recall that equivalence classes of star-products on $(M,\om)$
are parametrized by elements in
$$ \frac{1}{i \hb}  [\om] + H^2 (M, \C )[[\hb]] $$
called Fedosov characteristic classes. The following theorem was
proved by Karabegov and Schlichenmaier in \cite{ka} and \cite{KaSc},
in the case $K$ is the trivial line bundle. 
\begin{theo} \label{Fedosov}
The Fedosov class of the star-product $*_{\contravariant}$ is  $\frac{1}{i \hb} ( [\om] + \hb [\om_1])$.
\end{theo} 
Again it is interesting to note the appearance of $\om + \hb
\om_1$. We do not need this result but some related facts.
Let us define the {\em normalized symbol} of a Toeplitz operator by 
$$ \si_{\normalised} (T_k ) := (\identite + \tfrac{\hbar}{2} \Delta ) \si_{\contravariant} (T_k) $$
where $\Delta$ is the holomorphic Laplacian acting on $\Ci(M)$.
Actually we are only interested in the leading and second order terms
of $\si_{\normalised} (T_k)$ and modifying the definition of
$\si_{\normalised} (T_k)$ by a $O(\hb^2)$ term wouldn't change the statements
of our results.  To compare with our previous article \cite{oim2}, the Weyl symbol that
we introduced when $K$ is
the trivial line bundle is equal to the normalized symbol modulo $O(\hb^2)$.

The map $\si_{\normalised}: \Toeplitz \rightarrow \Ci(M)[[\hb]]$ satisfies
the same properties as $\si_{\contravariant}$ stated in theorem \ref{sc_algebra}.  Denote
by $*_{\normalised}$ the associated star product.

\begin{theo} \label{subprincipal_calculus} 
Let $f$ and $g$ belong to $\Ci(M)[[\hb]]$, then
$$ f *_{\normalised}  g = f.g + \tfrac{\hb}{2i} \langle \pi, df \wedge
dg \rangle + O(\hb^2) $$
and 
$$  i \hb^{-1} \bigl( f *_{\normalised}  g -  g *_{\normalised}  f
\bigr) =  \langle \pi + \hb \pi_1 , df \wedge
dg \rangle + O(\hb^2).$$
where $\pi$ is the Poisson bivector and $\pi_1$ is the bivector such that 
$\langle \pi_1, df \wedge dg
\rangle + \langle X_f \wedge X_g , \om_1 \rangle =0$
for every $f, g \in \Ci(M)$. 
\end{theo}

So $*_{\normalised}$ is a normalized star-product, in the sense that
the second order term in the first formula is antisymmetric, which explains
our terminology. 
Observe that $\pi + \hb \pi_1$ is the Poisson bivector associated to $\om
+ \hb \om_1$ in the sense that 
$$ \langle \pi + \hb  \pi_1, df \wedge dg
\rangle = \langle (X_f + \hb X_f^1) \wedge (X_g + \hb X_g^1) , \om +
\hb \om_1
\rangle +O(\hb^2),$$ 
where $X_f + \hb X_f^1+O(\hb^2)$  is
the Hamiltonian vector field of $f$  with respect to $\om +
\hb \om_1$, that is 
$$df + \langle \om + \hb \om_1 , X_f + \hb X_f^1
\rangle = O(\hb^2),$$ and the same holds for $g$ and $X_g + \hb
X_g^1$. So it follows from theorem \ref{Fedosov} that there exists a star
product equivalent to $*$ satisfying the formulas of theorem
\ref{subprincipal_calculus}. This last result is more precise because the
equivalence is specified. 
We can prove it  using the methods of \cite{oim1} or 
\cite{KaSc}. But this leads to complicated computations. We will present in
\cite{new} a more conceptual proof. We stated this result because we
can deduce from it a part of the Bohr-Sommerfeld conditions 
(cf. sections \ref{trace}, \ref{func}, \ref{var}). 

\scratch{Last remark is that $\si_{\normalised}$ isn't the only symbol map which
defines a normalized star-product. Indeed consider a symbol
map $\si$ of the form
$$ \si (T_k) = (\identite + \hb P ) \si_{\normalised} (T_k)  +O(\hb^2)$$
where $P$ is a differential operator. 
Then the associated star-product $*$ is normalized if and only if $P$ acts
as a vector field $X$.}

\subsection{Relation with geometric Quantization}

Our definition of the normalized symbol agrees in some sense with the
usual procedure to quantize observables in geometric
quantization. Assume that the Hamiltonian flow of $f \in \Ci (M)$
preserves the complex structure of $M$. Assume also that $K = L_1
\otimes \delta$, where $(\delta, \varphi)$ is a half-form bundle. 
Then the following operator is well-defined  
\begin{gather}
\label{geometric_quantization}
 \operateur (f) := f + \frac{1}{ik} \bigl( \nabla_{X_f}^{L^k \otimes L_1}\otimes
\identite + \identite \otimes \Lie_{X_f} \bigr) : \Hilbert_k \rightarrow
\Hilbert_k. \end{gather}
Here $\Lie_{X}$ acts on sections of
$\delta$ by $ \varphi ( (\Lie _X s) \otimes s) = \frac{1}{2} \Lie _X 
\varphi ( s^{\otimes 2} )$, or
equivalently as a Lie derivative where the pull-back of sections of
$\delta$ by a
complex diffeomorphism $\zeta $ is defined in such a way that $\varphi
((\zeta^* s)^{\otimes 2} )  =
  \zeta^* \varphi (s^{\otimes 2} )$.

The definition \eqref{geometric_quantization}
is natural for the following reason. Denote by $\Phi_t$ the
Hamiltonian flow of $X_f$. Let $\tilde{\Phi}_t$ be the lift
of $\Phi_t$ to $L^k \otimes L_1 \otimes \delta$
defined by the tensor product of the parallel transport along the trajectories of $X_f$ in $L^k
\otimes L_1$ and by the pull-back in $\delta$. Then the solution of the
Schr{\"o}dinger equation 
$$ \frac{1}{ik} \frac{d}{dt} \Psi(.,t) +  \operateur (f) \Psi (.,t) = 0$$
with initial condition $\Psi \in \Hilbert _k$ is given by 
$$\Psi (x,t) = e^{\frac{k}{i} \int_0^t f ( \Phi_{s-t} (x) ) \sl{ds}} \tilde{\Phi}_{t} \bigr( \Psi
(\Phi_{-t}(x)) \bigl) .$$

\scratch{La formule est v{\'e}rifi{\'e}e. Regarder dans le Woodhouse
s'il introduit qqch de ressemblant.} 

The important point for us is that $\operateur (f)$ is a Toeplitz
operators whose normalized symbol is $f$ modulo $O(\hb^2)$. We will prove a
more general result for every smooth function $f$, which simplifies
some further proofs. 

If $f$ is an arbitrary smooth function, formula \eqref{geometric_quantization} doesn't
make sense, because the Lie derivative with respect to $X_f$ doesn't
necessarily preserve $\Om^{n,0} (M)$. So we define for any vector field $X$ the
operator $D_X$, 
$$ D_X \al  = p (\Lie _X \al), \qquad  \al \in \Om^{n,0} (M) $$
where $p$ is the projection from $\La^{n} T^*M \otimes \C$ onto
$\La^{n,0} T^*M$ with kernel the sum
$$    \La^{n-1,1}T^* M \oplus \La^{n-2,2}T^* M \oplus ... \oplus
\La^{0,n} T^* M .$$ 
Next we let $D_X$ act on the sections of $\delta$, as the first-order
differential operator such that $ 2 \varphi( s \otimes D_X s) = D_X
\varphi( s^2)$ . 
 
\begin{theo} \label{normalised_symbol_quantisation}
For any $f\in \Ci(M)$, the operator 
$$  \operateur (f) := \Pi_k \Bigl( f + \frac{1}{ik} \bigl( \nabla_{X_f}^{L^k \otimes L_1}\otimes
\identite + \identite \otimes D_{X_f}  \bigr) \Bigr)  : \Hilbert_k \rightarrow
\Hilbert_k  $$
is a Toeplitz operator with principal symbol $f$ and vanishing subprincipal symbol. 
\end{theo}

This theorem is a consequence of the following lemma and an argument
of Tuynman \cite{Tu}. 

\begin{lemme} \label{l1}
Let $s$ be a half-form, then 
$$ D_{X_f} s = \nabla^{\delta}_{X_f} s + \tfrac{i}{2} (\Delta f) s $$
where $\nabla^\delta$ is the Chern connection of $\delta$ and $\Delta$
is the holomorphic Laplacian of $M$. 
\end{lemme}

It follows that 
\begin{gather} \label{defe}
  \nabla_{X_f}^{ L^k \otimes L_1}\otimes
\identite + \identite \otimes D_{X_f} = \nabla_{X_f}^{  L^k \otimes
  L_1 \otimes \delta } +  \tfrac{i}{2} (\Delta f) = \nabla_{X_f}^{  L^k \otimes
  K} +  \tfrac{i}{2} (\Delta f) \end{gather}
 Now we have for every $\Psi \in \Hilbert_k$, 
$$ \Pi_k ( \nabla_{X_f}^{  L^k \otimes
  K} \Psi)  = \tfrac{1}{i} \Pi_k ( \Delta f . \Psi ), $$
cf. \cite{Tu} or \cite{BoMeSc} for a proof. Hence
$$ \operateur (f) \Psi =  \Pi_k \bigl(
( f - \tfrac{1}{2k} \Delta f ) \Psi \bigr)
$$
which proves theorem \ref{normalised_symbol_quantisation}.

The last expression in \eqref{defe} shows that the definition of
$\operateur (f)$ is
independent of the choice of the half-form bundle and generalizes in
the cases where no such bundle exists.

\begin{proof}[Proof of lemma \ref{l1}]
It suffices to prove that for every $\alpha \in \Om^{n,0} (M)$, we have 
$$  D_{X_f} \al = \nabla_{X_f} \al + i (\Delta f) \al $$
Introduce normal complex coordinates $z^1,...,z^n$ centered at
$x_0$. So if $\alpha = dz^1 \wedge ... \wedge dz^n$, then $\nabla \al =
0$ at $x_0$. Let us write $\om = i G_{j,k} dz^j \wedge d
\bar{z}^k$. Then 
$$X_f = - i G^{j,k} (\partial_{z^j} f) \partial_{\bar{z}^k} + i
G^{j,k}(\partial_{\bar{z}^k} f) \partial_{z^j}
$$ 
Using that the first derivatives of $G^{j,k}$ vanish at $x_0$, we
obtain easily 
$$ D_{X_f} \al = i G^{j,k}  (\partial_{z^j}\partial_{\bar{z}^k} f) = i
\Delta f.$$
The result follows. 
\end{proof}

\section{Bohr-Sommerfeld Conditions} \label{BS_section}

\subsection{The result}

Assume that $M$ is $2$-dimensional. Let $(\delta,\varphi)$ be a
half-form bundle and let us write $K = L_1 \otimes \delta$ as
previously. Consider a self-adjoint Toeplitz
operator $(T_k)$. Its normalized symbol $$f_0 + \hb f_1 +...$$  is
real-valued. Bohr-Sommerfeld conditions give the spectrum of $T_k$ on
every open interval $I$ of regular values of $f_0$ in the
semi-classical limit. To simplify the statements, assume that $f_1$ vanishes. 

Let $\Ga^1$, ...,
$\Ga^m$ be the components of $f_0^{-1}(I)$. For every $i \in
\{ 1,...,m\}$, the map $$f_0 : \Ga^i \rightarrow I$$ is a trivial
fibration with fiber diffeomorphic to $
S^1$. For every $\la \in I$, fix an orientation on the fiber
$\Ga^{i}_\la = f_0^{-1} (\la) \cap \Ga^i$ depending continuously on
$\la$.

Let $a^i \in \Ci(I)$ be the {\em principal action}, defined in such a
way that the parallel transport in $L$ along $\Ga^i_\la$ is the
multiplication by $\exp \bigl( {i a^i(\la)}\bigr)$. Using $L_1$
instead of $L$, define in the same way the
  {\em subprincipal action} $a^i_1 \in \Ci(I)$.

Let us define an index $\epsilon^i$ from the half-form bundle
$(\delta, \varphi)$. Observe that the restriction $\delta_{i, \la}$ of
$\delta$ to $ \Ga^{i}_\la$ is a square root of $T^*
\Ga^{i}_\la \otimes \C$. Indeed, let us denote by $\iota$ the embedding $\Ga^{i}_\la
\rightarrow M$, then the map
$$ \varphi_{i,\la} :\delta_{i, \la}^2 \rightarrow T^*
\Ga^{i}_\la \otimes \C, \qquad u \rightarrow \iota^* \varphi (u)$$
is an isomorphism of line bundle. The set 
$$ \{ u \in \delta_{i,\la}; \; \varphi_{i,\la} (u ^{\otimes 2}) >0 \} $$
has one or two connected components. In the first case,
we set $\epsilon^{i}_\la =1$ and in the second case $\epsilon^{i}_\la
=0$. Observe that $\epsilon^{i}_\la$ doesn't depend on $\la$.

The Bohr-Sommerfeld conditions are   
\begin{gather} \label{BS_conditions}
  a^i(\la) + k^{-1}( a^i_1 (\la) + \epsilon^i \pi   ) \in \frac{2 \pi}{k} \Z \end{gather} 
Denote by $\Si^i(k)$ the set of $\la \in I$ satisfying
\eqref{BS_conditions}. When $k$ is sufficiently large, $\Si^i(k)$ is a
finite set containing 
$$\frac{k}{2\pi} \volume (\Ga^i ) + O( k^{-1})$$ points. 
Let $\Si(k)$ be the union of the $\Si^i(k)$. Define the
{\em multiplicity} of $\la \in \Si(k)$ as the number of $\Si^i(k)$ which contains
$\la$. The points of $\Si (k)$ approximate the eigenvalues of $T_k$
in the following sense.

\begin{theo} \label{BS_theo}
 Let $\la_-(k)$ and $\la_+(k)$ be two sequences of $I$ such that 
\begin{gather} \label{E1}
d ( \la_-(k) , \Si (k)) \geqslant  C k^{-1}, \quad d ( \la_+(k) ,
\Si (k))
\geqslant C k^{-1} \end{gather} 
for some positive $C$. Assume furthermore that there exists $\la_-,
\la_+ \in I$ such that 
\begin{gather} \label{E2}
\la_- \leqslant
\la_-(k) \leqslant \la_+(k) \leqslant \la_+ .\end{gather} 
Denote by $\la_1(k) \leqslant \la_2(k) \leqslant ....\leqslant \la_{N(k)}(k)$
 (resp. $\la'_1(k) \leqslant \la'_2(k) \leqslant ....\leqslant \la'_{N'(k)}(k)$)
 the eigenvalues of $(T_k)$ (resp. points of $\Si (k)$) contained in
 $(\la_-(k), \la_+(k))$ and counted with multiplicities. Then, when $k$ is
 sufficiently large, $N(k) = N'(k)$. Furthermore  
\begin{gather} \label{Approximation} 
 \la_j(k) =  \la_j' (k) + O(k^{-2}) \end{gather}
uniformly with respect to $j$.  
\end{theo} 

The interest of condition \eqref{E1} is to avoid any ambiguity in the counting of eigenvalues near the
endpoints of $(\la_- (k), \la_+ (k))$. It is not restrictive. Indeed  if
$\la_-(k)$, $\la_+ (k)$ are arbitrary sequences satisfying \eqref{E2},
then by modifying them by suitably chosen $O(k^{-1})$ sequences, we
obtain sequences satisfying both estimates \eqref{E1} and \eqref{E2}. 

Since the definition of the Toeplitz operators and of their normalized
symbol only depend on $K$ and not on the choice of the half form bundle, it is likely that the same holds for the Bohr-Sommerfeld
conditions. This is easily checked using that any other half-form bundle is of the form $\delta' =\delta \otimes
F$, where $F$ is a flat Hermitian line bundle with holonomy in $\Z_2$. So
$L_1' = L_1 \otimes F^{-1}$, and straightforward computations show
that the functions $a^i_1 +
\epsilon^i$ do not depend on the choice of $\delta$. 

To compare with our previous results in \cite{oim2}, when $K$ is
the trivial bundle, we defined the function $a^i_1 +
\epsilon^i$ as the integral of the geodesic curvature of
$\Ga_{\la}^i$. 

We can also approximate the eigenvalues up to a $O(k^{-\infty})$ error. More precisely,  there exist
sequences $(S^i(.,k))_k$ of $\Ci (I)$ such that the Bohr-Sommerfeld
conditions $$S^i(.,k) \in \tfrac{2 \pi}{ k} \Z$$ 
instead of \eqref{BS_conditions} lead to the same
result with a $O(k^{-\infty})$ error in \eqref{Approximation}. Furthermore, the sequences $S^i(.,k)$ admit asymptotic
expansions of the form 
$$S_0^i + k^{-1} S^i_1 + k^{-2} S^i_2 +...$$ 
The leading and second order terms are the same like in \eqref{BS_conditions}. Applying an argument of Colin de
Verdi{\`e}re \cite{Co}, we can prove that the derivatives of the $S^i_j$ only depend on the star-product $*_\normalised$ and the
normalized symbol of $T_k$.  In particular, assuming that $S_0^i$ is
increasing, we will deduce from
proposition \ref{subprincipal_calculus} that necessarily
\begin{gather} \label{Colin}
S^i_0 (\la) -   S^i_0 (\la') = \int_{D^i} \om \quad \text{and} \quad S^i_1 (\la ) -  S^i_1 (\la') = \int_{D^i} \om_1 \end{gather}
if $\la , \la' \in I$ are such that $\la' \leqslant \la$ and $D^i =
\Ga^i \cap f_0^{-1} ((\la , \la'))$.
This determines the Bohr-Sommerfeld conditions \eqref{BS_conditions} modulo a
constant term. 
Note that only
the derivatives of $S^i_0$ and $S^i_1$ are determined by the
star-product. Indeed, if we twist $L$ or $L_1$ by a flat Hermitian
line bundle, the actions of the non-contractible loops may change although the star-product remains the same. 

As a first step to deduce \eqref{Colin} from proposition
\ref{subprincipal_calculus}, let us state some results on traces
and functional calculus of Toeplitz operators.

\subsection{Traces} \label{trace}

In this section and the next one, we do not necessarily assume that $M$
 is $2$-dimensional. 
 It is a known result that the trace of a Toeplitz operator $T_k$ with normalized symbol $f_0 +
\hb f_1 +...$ admits an asymptotic
expansion of the following form:
$$ 
 \Trace ( T_k) = \Bigl( \frac{k}{2 \pi} \Bigr)^n \int_M (f_0 +
k^{-1} f_1 +...) (1+ k^{-1} d_1 + k^{-2} d_2 +...) \mu_M
+O(k^{-\infty}) $$
where $d_1$, $d_2$, $d_3$, ... are functions of $\Ci(M)$ which do not depend of
$T_k$ (cf. for instance \cite{oim1}). 

These functions may be computed in terms of the K{\"a}hler metric
 by using the methods of \cite{oim1}, but it is more convenient to relate them to the star-product $*_\normalised$. 
To do this, observe that the $\C[[\hb]]$-linear map
\begin{gather} \label{Trace}  \trace: \Ci(M) [[\hb ]] \rightarrow \C[\hb^{-1}, \hb]], \quad f(\hb)
\rightarrow (2 \pi \hb )^{-n} \int_M f(\hb) d(\hb)\; \mu_M \end{gather}
where $d( \hb) =1 + \hb d_1 + \hb^2 d_2 +...$ is a trace for the
star-product $*_{\normalised}$, in the sense that it satisfies 
$$\trace (f
*_{\normalised} g) = \trace( g *_{\normalised} f).$$ 
Following Fedosov
\cite{Fe} or Nest-Tsygan \cite{NeTs}, such a trace is
unique up to multiplication by an element of $\C [\hb^{-1}, \hb]]$ and
there exists a canonical one determined by the
following normalization condition: for every local equivalence $\Phi$
between  $*_{\normalised}$ and the
Weyl star-product, we have 
\begin{gather}  \label{trace_normalised}
\trace (f) =  (2 \pi \hb)^{-n} \int_M \Ph (f)\; \mu_M, \quad \text{
  where } \mu_M = \om^n/n!. \end{gather} 
We claim that the trace defined in
\eqref{Trace} is the canonical trace. This follows from the fact that
the quantization by Toeplitz operators is microlocally equivalent to
the usual Weyl quantization.

So the functions $d_i$ are determined by $*_\normalised$. In
particular, it follows from
proposition \ref{subprincipal_calculus} that 
\begin{gather} \label{sub_trace} 
 \trace (f) =  (2 \pi \hb)^{-n} \int_M f \; ( \om + \hb \om_1 )^{\wedge
  n}/ n! + O( \hb ^{-n +2}) \end{gather}  

\begin{proof}[Proof of formula \eqref{sub_trace}] Consider an equivalence of star-product of the form
  $$\Phi = \identite + \hb X + O(\hb^2)$$ where $X$ is a vector field. It is easily
  checked that the star-product 
$$ f *' g =  \Phi \bigl( \Phi^{-1} (f) *_{\normalised} \Phi^{-1} (g) \bigr)$$
is normalized and satisfies
$$  i \hb^{-1} \bigl( f *'  g -  g *'  f
\bigr) =  \langle \pi + \hb( \pi_1+ X.\pi) , df \wedge
dg \rangle + O(\hb^2).$$
Locally on can choose $X$ in such a way that $ \pi_1 + X.\pi = 0$. Indeed this
equation is equivalent to 
$\om_1 + d\al = 0,$ where $\al$ is the
$1$-form such that $ \om(X,.) = \al$. 
 
Next step is to introduce local Darboux coordinates, which define a Weyl star-product $*_{\operatorname{Weyl}}$. And  modifying $\Phi$
by a $O(\hb^2)$ term, one has $*_{\operatorname{Weyl}} = *'$
(cf. \cite{BeCaGu}).  Then it follows from
\eqref{trace_normalised} that 
$$ \trace (f) = (2 \pi \hb)^{-n} \int_M f\; \bigl(\mu_M - \hb (X\mu_M)
+O(\hb^2)\bigr) .$$
By definition of $\pi_1$ in theorem
\ref{subprincipal_calculus}, one has  $\langle \pi_1 , \om \rangle +
\langle \pi, \om_1 \rangle = 0 $. Since $X. \langle \pi, \om \rangle
=0$, we obtain 
$$\langle \pi , X \om \rangle = -  \langle X  \pi, \om
\rangle = \langle \pi_1, \om \rangle = - \langle \om_1, \pi \rangle.$$
Then a straightforward computation leads to 
$$\langle \pi^ {\wedge n}, \mu_M - \hb X
\mu_M \rangle  = \langle \pi^{\wedge n}, ( \om + \hb \om_1 )^{\wedge
  n}/ n! \rangle + O(\hb^2),$$ 
which proves the result.
\end{proof}

As a consequence of \eqref{sub_trace}, we obtain the estimate
\eqref{Riemann-Roch} of the dimension of $\Hilbert_k$, since this
dimension is the trace of the projector $\Pi_k$, whose normalized symbol is
$1$.
 Actually the index theorems of deformation quantization proved
in \cite{Fe} or \cite{NeTs} yield the asymptotic expansion of $\Trace (
\Pi_k)$ modulo $O(k^{-\infty})$ in terms of the Fedosov class
of $*_\normalised$. Since this trace is an integer, the
$O(k^{-\infty})$ error vanishes when $k$ is sufficiently
large. In this way we can deduce Riemann-Roch-Hirzebruch theorem from
theorem \ref{Fedosov} .

\subsection{Functional calculus} \label{func}

Let $T_k$ be a self-adjoint Toeplitz operator with normalized symbol
$$f =f_0 + \hb f_1 +....$$ and let $g$ be a function of $\Ci
(\R, \C)$. Then it is known that $g(T_k)$ is a Toeplitz
operator (cf. for instance \cite{oim1}). Furthermore, the normalized symbol of $g(T_k)$ is given by
the following non-commutative Taylor formula: 
\begin{gather} \label{NC_Taylor}
 g^{*_{\normalised}}(f) (x)= \sum \frac{1}{\ell!}g^{(\ell)} (f_0 (x)) \; ( f(y) - f_0 (x)
)^{*_{\normalised} \ell} \eval{y = x}, \end{gather} 
where $g^{(\ell)}$ is the $\ell$th derivative of $g$ and $h^{*_{\normalised}
  \ell} = h*_{\normalised}...*_{\normalised} h$ repeated $\ell$ times. 
In particular, an easy computation
from proposition \ref{subprincipal_calculus} leads to
\begin{gather} \label{sub_functional_calculus} 
g^{*_{\normalised}}(f)=  g(f_0) + \hb g'(f_0)f_1 +
O(\hb^2). \end{gather} 
\scratch{
Notice that formula \eqref{NC_Taylor} makes sense for an arbitrary
star product $*$. This defines a formal functional calculus
satisfying the expected properties. For instance, one can check that $g_1^* (f) * g_2^*(f) = (g_1
g_2)^{*}(f)$. }

\subsection{On the variation of $S_0$ and $S_1$} \label{var}
 
Let us deduce formulas \eqref{Colin} on the variations of $S_0$ and
$S_1$ from the trace formula \eqref{sub_trace} and the functional symbolic calculus. 
Assume that
$f_0^{-1}(I)$ is connected. 

If  $g \in \Ci_o (I, \C)$, then we deduce from
\eqref{sub_functional_calculus}, \eqref{sub_trace} and the fact that
$f_1 =0$ that 
\begin{xalignat}{1} \notag 
\Trace (g(T_k)) = & \frac{k}{2 \pi} \int_M (g \circ f_0)  ( \om + k^{-1}
\om_1 ) \; + O(  k^{-1})\\\label{eq1}
=& \frac{k}{2 \pi} \int_I g \; \bigl( f_{0*}  [ \om + k^{-1}
\om_1 ] \bigr) \; + O( k^{-1})
 .\end{xalignat}
where  $f_{0*}$ is the push-forward $\Om^2 (M) \rightarrow \Om^1 (I)$
defined by 
$$  \int_M (f_0^* h ) \ \al  = \int_I g. f_{0*} \al, \qquad \forall h
\in \Ci_o(I) .$$

On the other hand, assume the Bohr-Sommerfeld condition is $$ S(\la,k) \in 2 \pi k^{-1} \Z $$ with
$S(\la,k) =S_0 + k^{-1}S_1 + O(k^{-2})$. If the derivative of $S_0$
doesn't vanish, then one can invert the functions $S(.,k)$ when $k$ is
sufficiently large and
$$ \Trace (g(T_k)) = \sum_{x \in 2\pi k^{-1} \Z} g ( S^{-1}(x,k))
+ O(k^{-1}) $$
Interpreting this as a Riemann sum, it follows that
\begin{xalignat}{1} \notag
 \Trace (g(T_k)) =& \frac{k}{2 \pi} \int_\R g (S^{-1}(x,k )) dx
 +O(k^{-1}) \\ \label{eq2}
 = &  \frac{k}{2 \pi} \int_I g (\la) S' (\la, k ) d\la
 +O(k^{-1})
\end{xalignat}
with the same orientation of $I$ as before if $S'_0$ is positive.  
Since \eqref{eq1} and \eqref{eq2} hold for any function $g$ of $\Ci_o( I,
\C)$, we have 
$$ d S_0   + \hb d S_1  =  f_{0*}  ( \om + \hb \om_1 )$$ 
and \eqref{Colin} follows.

\section{Lagrangian states}

\subsection{Definitions and symbolic calculus} \label{sec:resultat}

First we recall the definition of a local Lagrangian section
associated to a closed Lagrangian embedding $\iota: \Ga \rightarrow M$. 

Let $U$ be an
open set of $M$ such that $U_\Gamma := \iota^{-1} (U )$ is contractible. Since the
curvature of $\iota^* L$ vanishes, there exists a 
flat unitary section $t_\Gamma$ of $ \iota^* L \rightarrow U_\Gamma$. 
Introduce a formal
series 
$$\sum_{\ell =0 }^{\infty} \hb^\ell g_\ell \in \Ci ( U_\Gamma, \iota^*K)[[\hb]].$$ 
Let $V$  be an open set of $M$ such that $\overline{V} \subset U$. 
Then a
 sequence $\Psi_k \in {\mathcal{H}}_k$ is a {\em Lagrangian
 section}
 over $V$ associated to $(\Gamma, t_\Gamma)$ with symbol $\sum \hb^\ell
 g_\ell$ if 
\begin{gather*} 
 \Psi_k(x) =  \Bigl( \frac{k}{2\pi} \Bigr)^{\frac{n}{4}} F^k (x)
 \tilde{g} (x,k) + O(k^{-\infty}) \text{ over $V$, } \end{gather*} 
where 
\begin{itemize} 
 \item  $F$ is a section of $L \rightarrow U$ such that 
$$\iota^* F =
 t_\Gamma \quad \text{ and } \quad \bar{\partial} F \equiv 0$$ modulo a section which
 vanishes to every order along $\iota (\Gamma)$. Furthermore $|F(x)|< 1$ if $x \notin \iota(\Gamma)$.  
 \item
$\tilde{g}(.,k )$ is a sequence of $\Ci ( U,K)$
with an asymptotic expansion $\textstyle{\sum} k^{-\ell} \tilde{g}_{\ell}$
in the $\Ci$ topology such that 
$$\iota^* \tilde{g}_\ell = g_\ell \text{
 and } \quad \bar{\partial} \tilde{g}_\ell \equiv 0$$ modulo a section which vanishes at
every order along $\iota(\Gamma)$.
\end{itemize}
We assume furthermore that $\Psi_k$ is admissible in the sense that $\Psi_k(x)$ is uniformly
$O(k^{N})$ for some $N$ and the same holds for its successive
covariant derivatives.

It is not obvious that such a sequence exists.
\begin{theo} \label{existence}
For every series $\sum \hb^\ell g_\ell$ of $\Ci ( U_\Gamma,
 \iota^* K)[[\hb]]$, there exists a Lagrangian
  section  over $V$ associated to $(\Gamma, t_\Gamma)$ with symbol $\sum \hb^\ell
 g_\ell$. It is unique modulo a section which is $O(k^{-\infty})$
  over $V$. \end{theo}

In the statement of the following theorems, we consider that $K = L_1
\otimes \delta $ over $U$, where $(\delta, \varphi)$ is a half-form bundle. Recall that
$\iota^* \delta = \delta_{\Ga}$ is a square root of $\La ^n T^*\Ga \otimes
\C$ through the isomorphism 
$$ \varphi_{\Ga} : \delta^2_ \Ga
\rightarrow \La ^n T^*\Ga \otimes \C, \quad u \rightarrow \iota^*
\varphi( u).$$ 
Let us associate
to the {\em principal symbol} $g_0$ of a Lagrangian section a density
$m(g_0)$ where $m$ is the map:
$$ m : \iota^* L_1 \otimes  \delta_{\Ga} \rightarrow | \La | (\Ga) , \quad
u \otimes v \rightarrow \| u \|^2_{L_1} \ | \varphi_{\Ga} ( v^{\otimes
  2})|$$
We then have the following estimate of the norm of $\Psi_k$. 
\begin{theo} \label{norme} 
Let  $\xi \in \Ci_o ( V)$, then we have  
$$ \int_M \xi \  \| \Psi_k \|^2 _{L^k \otimes L_1 \otimes \delta } \; \mu_M
 = \int_{\Gamma} (\iota^* \xi) \  m(g_0) +
 O(k^{-1}).$$
\end{theo}
Next results describe how a Toeplitz operator acts on a Lagrangian
section.

\begin{theo}  \label{Toep_Lag}
Let $T_k$ be a Toeplitz operator with principal symbol
$f_0$. Then $T_k \Psi_k$ is a Lagrangian section over $V$ associated
to $(\Gamma, t_\Gamma)$ with symbol $(\iota^* f_0) g_0 +O(\hb)$. 
\end{theo}

To prove the Bohr-Sommerfeld conditions, we need to compute
the subsequent coefficient of the symbol of $T_k \Psi_k$, in the case
where $f_0$ is constant over $\Ga$. 
 
\begin{theo} \label{resultat_clef} 
Let $T_k$ be a Toeplitz operator with normalized symbol $f_0 + \hb
 f_1 + O(\hb^2)$. Assume that $f_0$ is constant along $\Gamma$. Then
the symbol of $ T_k \Psi_k$ is 
$$ (\iota^* ( f_0 + \hb f_1)).( g_0  + \hb g_1)  + \hb  \tfrac{1}{i} (
\nabla_{X}^{\iota^* L_1}
\otimes \identite + \identite \otimes \Lie^{\delta_{\Ga}}_X ).g_0  + O(\hb^2) $$
where 
\begin{itemize} 

\item $X$ is the Hamiltonian vector field of $f_0$,
\item  $\nabla^{\iota^* L_1}$ is the pull-back of the Chern connection of $L_1$, 
\item $ \Lie^{\delta_{\Ga}}_X$ is the first order differential operator
  acting on sections of $\delta_{\Ga}$ 
  such that $$ \varphi_{\Ga} ( \Lie^{\delta_{\Ga}}_X g  \otimes g)  = \tfrac{1}{2}
\Lie_X \varphi_{\Ga} ( g ^{\otimes 2}) $$
for every section $g$.
\end{itemize} 
\end{theo}

It is easily checked that the operator $\nabla_{X}^{\iota^* L_1}
\otimes \identite + \identite \otimes \Lie^{\delta_{\Ga}}_X$  doesn't
depend of the choice of the half-form bundle if we consider that it
acts on sections of $\iota^* K = \iota^* L_1 \otimes
\delta_\Ga$. The
same holds with the map $m$. 

\subsection{Proof of Bohr-Sommerfeld conditions}
Let us deduce from the previous theorems the Bohr-Sommerfeld conditions for $n$
self-adjoint commuting Toeplitz operators $T^1$, $T^2$,..., $T^n$, which is  a slight
generalization of \eqref{BS_conditions}. 

Denote by $f^i_0$ and $f^i_1$
the principal and subprincipal symbols of $T^i$. Let $E$ be a
regular value of $f=(f^1_0,...,f^n_0)$ and $\iota : \Ga \rightarrow M$ be
an embedding with image a connected component
of $f^{-1} (\{ E\})$. Since $M$ is compact, $\iota (\Ga)$ is a Lagrangian
torus. So there exists a half-form bundle $(\delta, \varphi)$ defined over
a neighborhood of $\iota (\Ga)$. It is not unique but as usual the
final result doesn't depend on the choice of $(\delta, \varphi)$. 
Introduce like in the previous section two open sets $U,V$ and a flat section
$t_\Ga$. Let us try to solve the eigenvalues equation 
\begin{gather} \label{vap_equation} 
 T^i \Psi = E^i \Psi + O(k^{-\infty}) \text{ over $V$} \end{gather}  
where $\Psi$ is a Lagrangian section associated to $(\Ga, t_{\Ga})$. By
theorem \ref{resultat_clef}, the symbol of $(T^i
-E^i) \Psi$ is $O(\hb)$ because $\iota ^* f^i = E^i$. Furthermore it
is $O(\hb^2)$ if and only if it satisfies the following transport equation
\begin{gather} \label{transport_equation}
 \bigr[  f^i_1 + ( \nabla_{X^i}^{\iota^*L_1  }
\otimes \identite + \identite \otimes \Lie^{\delta_{\Ga}}_{X^i} ) \bigl] g_0 = 0
\text{ over $V \cap \Ga$} \end{gather}
where $g_0$ is the principal symbol of $\Psi$. This equation can be
interpreted as $g_0$ being flat for a connection on $ \iota^* L_1 \otimes \delta_{\Ga}$ that we
describe  now. 

First a section $g$ of $\delta_{\Ga}$ is flat if $\mathcal{L}_{X^i}^{\delta_{\Ga}} g
= 0$ for every $i$. With this definition, $\varphi_{\Ga} : \delta_{\Ga}^2
\rightarrow \La ^n T^*\Ga \otimes \C $ is a morphism of flat bundles, if
we endow $\La ^n T^* \Ga \otimes \C $ with the Weinstein connection. Since
the form $\be$ such that 
$$ \langle  \be,   X^1 \wedge ... \wedge X^n \rangle =1 $$
is a global non vanishing flat section of $\La ^n T^* \Ga \otimes \C $,  $\delta_\Ga$ has  holonomy in $\Z_2 \subset U(1)$.

Let $\al \in \Om^1(\Gamma)$ be such that $ \frac{1}{i} \langle \al ,
X^j \rangle =  f^j_1$. Consider the connection  $\nabla ^{ \iota^* L_1} + \al$
on $\iota^* L_1$. Its flat sections satisfy 
$$(f_1^j + \tfrac{1}{i} \nabla^{\iota^* L_1}_{X^j}) s = 0.$$ 
Furthermore its curvature vanishes. Indeed, since $[X^i, X^j] = 0$, we have 
\begin{xalignat*}{2}
 \frac{1}{i}  \langle d \al , X^i \wedge X^j \rangle = &  ( X^i. f_1^j - X^j. f_1^i
  ) \\
= &  \om_1 ( X^i, X^j) \end{xalignat*}
which follows from proposition \ref{subprincipal_calculus} and the fact that
$[T^i, T^j] = 0$. 
 
This defines a structure of flat line bundle for $\iota^*L_1 \otimes
\delta_{\Ga} $, whose flat sections are the solutions of \eqref{transport_equation}.
Recall that 
$$ \Psi(x,k) =  t^k_\Ga (x)  (g_0(x) + O(k^{-1}))
\text{ over $V \cap \iota(\Gamma)$},$$
where $t_\Ga$ is a flat section of $\iota^* L $. The condition
to patch together these sections along $\Ga$ is the Bohr-Sommerfeld
condition:
$$ \text{$\iota^*( L^k \otimes   L_1) \otimes \delta_{\Ga} \rightarrow \Ga$ is trivial as a flat bundle.} $$
When $M$ is two-dimensional, this is equivalent to
\eqref{BS_conditions}. To prove theorem \ref{BS_theo} or a similar result in the
$2n$-dimensional case for the joint spectrum, we should consider
Lagrangian sections depending continuously on $\Ga$. Furthermore, we
can show by using a local normal form that the solutions of
\eqref{vap_equation} are necessarily Lagrangian sections associated to
$\Ga$. A complete proof is in \cite{oim2}. The only novelty here is
the formulation of theorems \ref{resultat_clef} and \ref{norme}, and
consequently of the Bohr-Sommerfeld conditions. 

\subsection{Comparison with the cotangent case} 
To compare theorems \ref{norme} and \ref{resultat_clef} with the similar statements in
  the case of pseudo-differential operators, we can introduce some
  kind of Maslov
  bundle in the following way. Recall that we denote by $\varphi_{\Ga}$ the
  isomorphism $\delta_{\Ga}^2 \rightarrow \La ^n T^*\Ga \otimes
  \C$. Introduce
$$ P := \bigl\{ u \in \delta_{\Ga} ; \; \varphi_{\Ga} ( u ^{\otimes 2} )
\in \La^n T^* \Ga - \{0 \} \bigr\} .$$
Let $\Z_4$ be the subgroup $\{ 1, -1, i, -i \}$ of $\C^*$. Then $P$ is
a principal bundle with structure group $ \Z_4 \times \R_+$. Introduce
the complex line bundles $| \delta_{\Ga} |$ and $\arg ( \delta_{\Ga}
)$ associated to $P$ via the homomorphism $\Z_4 \times \R_+ \rightarrow
\Z_4$ and $\Z_4 \times \R_+ \rightarrow \R_+$ respectively. Following Weinstein in \cite{We}, we call $\arg (\delta_{\Ga})$ the unitarization
of $\delta_{\Ga}$. We have a
canonical isomorphism
$$ \delta_{\Ga} \rightarrow | \delta_{\Ga} | \otimes \arg
(\delta_{\Ga} ).$$ 
Furthermore, the map
$$ |\delta_{\Ga} | = P
\times_{\R_+} \C \ni  [u,z] \rightarrow z.| \varphi_{\Ga} (u^{\otimes
  2})|^{\frac{1}{2}}\in  |\La |^{\frac{1}{2}} ( \Ga)$$
is an isomorphism between $| \delta_{\Ga} |$ and
the bundle of half-densities of $\Ga$. So we obtain an isomorphism
$$ \zeta: \delta_{\Ga} \rightarrow |\La |^{\frac{1}{2}} ( \Ga) \otimes \arg
(\delta_{\Ga} )  .$$
The bundle $\arg (\delta_{\Ga})$ is a line bundle with structure group $\Z_4$ like the
Maslov bundle. 
The isomorphism $\zeta$ intertwines the operator
$\Lie_X^{\delta_{\Ga}}$ of theorem \ref{resultat_clef} with the Lie
derivative of half-densities. So in the case $L_1$ is trivial, the
theorem \ref{resultat_clef} is similar to the formula 1.3.13 in
\cite{Du3} p223, computing the symbol of an oscillatory integral acted
on under a differential operator. 
Furthermore the map $ \delta_{\Ga}
\rightarrow |\La |(\Ga)$ used in  theorem \ref{norme} is the composition of
$\zeta$ with the squaring map from half-densities into
densities. Again theorem \ref{norme} is similar to formula 1.3.15 in
\cite{Du3} p224, computing the norm of an oscillatory integral.

To end this comparison, we apply the previous construction to a
symplectic vector space $E$ and prove that we obtain the usual Maslov
bundle. Consider a one-dimensional vector space $\delta$ with an
isomorphism $\varphi: \delta^2 \rightarrow \La^{n, 0} E^*$. 
Let $\Lag (E)$ be the Lagrangian Grassmannian and $\eta
\rightarrow \Lag (E) $ be the tautological vector bundle, that is the
bundle whose fiber over $x$ is $x$ itself. In the same way we defined
$\varphi_{\Ga}$, one has an isomorphism
$$ \varphi_{\Lag} : \Lag (E) \times \delta^2 \rightarrow \La ^n \eta
\otimes \C,  $$ 
sending $(x,u)$ into $\iota_x^* \varphi(u)$ where $\iota_x$ is the
embedding $ \eta_x \rightarrow E$. Introduce the
$\Z_4 \times \R_+$ bundle 
$$  \bigl\{ (x,u) \in \Lag (E)
\times \delta; \;  \varphi_{\Lag} (x,u^2) \in \La_x^n \eta \text{ and }
u \neq 0 \bigr\} $$
Dividing by $\R_+$ we get a $\Z_4$-principal bundle $\Mas$. We claim
that the holonomy in
$\Mas$ of a loop of $\Lag (E)$ is the mod 4 reduction of its Maslov index.

\begin{proof}  Introduce linear
Darboux coordinates $(p^i,q^i)$ and identify $E$ with $\R^{2n}$. Set $z^j= p^j + i q^j$ and let 
$$ \delta := (dz^1\wedge...\wedge dz^n)^{\frac{1}{2}} \C$$ 
 be the square root of $\La^{n,0} E = (dz^1\wedge...\wedge dz^n)\
\C$. Recall that $\Lag (E)$ is isomorphic to $U(n)/O(n)$ (cf. lemma 2.31 of
\cite{McSa}), through the map sending the unitary matrix $U =P +i Q$
 into the range of  
\begin{gather*}  
 A_U = \left( 
\begin{array}{c} P  \\ Q 
\end{array} \right) \end{gather*}  
Let us denote by $\al^1_U$,...,$\al^n_U$ the base of
$\eta_x^*$ dual to the column vectors of $A_U$. A
straightforward computation shows that $\varphi_{\Lag}$ sends the
square of 
$$([U], (dz^1\wedge...\wedge dz^n)^{\frac{1}{2}})
\in \Lag (E) \times \delta $$ into
$\det (U) \ \al^1_U \wedge... \wedge \al^n_U$. Consequently, 
$$ \Mas \simeq \bigl\{ ( [U], v ) \in \Lag (E)
\times \C^* ;
\; U \in U(n) \text{ and } v^2 \det ( U) = \pm 1 \bigr\} $$
Recall now that the
Maslov index  of a loop  $x: \R/ \Z \rightarrow \Lag (E)$ is the degree of $\rho \circ x : S^1 \rightarrow
S^1$ where $\rho$ is the map 
$$ \Lag (E) \rightarrow S^1, \quad [U] \rightarrow \operatorname{det}^2 (U)$$ 
(cf. \cite{McSa} page 53). Its mod 4 reduction is the holonomy of $x$
in $ \Mas$.
\end{proof}

This result is related to the paper
\cite{We} of Weinstein, where it is observed that the Maslov bundle of
$\Lag (E)$ is a unitarization of a square root of $\La ^n \eta
\otimes \C$.

Last remark is that in general the
Maslov bundle of a Lagrangian submanifold of a cotangent space can be
different of the bundle we construct. Indeed notice that the structure
group of $\arg (\delta_{\Ga})$ reduces to $\Z_2$ if and only if $\Ga$
is orientable. Consider a non-orientable manifold $Q$. Then as the
null section of $T^*Q$, $Q$ is a
Lagrangian submanifold and its Maslov bundle is the flat trivial
bundle. So it can not be a unitarization of a square root of $\La^n T^*Q
\otimes \C$.

\section{Proof}

We assume in the whole section 
that there exists a globally defined half-form bundle $(\delta, \varphi)$ and
$K = \delta$. There is no difficulty to generalize to the case where
$K = \delta \otimes L_1$. 

\subsection{A preliminary result} 

Consider a sequence $\Psi_k \in \Ci(M,L^k \otimes \delta)$ of the
form 
$$\Psi_k(x) =  \Bigl( \frac{k}{2\pi} \Bigr)^{\frac{n}{4}} F^k (x)
 \tilde{g} (x,k) + O(k^{-\infty}) \text{ over $V$, } $$
where $F$ and $\tilde{g}(.,k)$ 
  satisfies the same assumptions as in section \ref{sec:resultat} except that the coefficients
  $\tilde{g}_k$ do not necessarily satisfy $\bar{\partial} \tilde{g}_k \equiv
  0$.  Assume furthermore that $\Psi_k$ is admissible. 

\begin{theo} \label{theo:main}
Let $T_k$ be a Toeplitz operator with principal symbol
  $f_0$. Then $T_k \Psi_k$ is a Lagrangian section over $V$ with
  symbol $\iota^* (f_0 \tilde{g}_0) + O(\hb)$. 
 
Furthermore, if $\tilde{g}_0$
  and its first derivatives vanish along $\iota(\Ga)$, then the
  symbol of $T_k \Psi_k$ is
$$   \hb (\iota^* f_0) \bigl(  \square \tilde{g}_0 + \iota^* \tilde{g}_1) ) + O(\hb^2)
$$
where $\square \tilde{g}_0 \in \Ci(\Ga, \delta_{\Ga})$ and at every $x \in \Ga$ 
$$  \square \tilde{g}_0  (x) = - \frac{1}{2}  \sum \bar{\partial}_i \bar{\partial}_i
\tilde{g}_0 (\iota(x))$$
if $\partial_1,...,\partial_n$ is a base of vectors of
$T^{1,0}_{\iota(x)} M$ such that  
$\frac{1}{i} \om (\partial_{i}, \bar{\partial}_j)
= \delta_{ij}$
and the vectors $\partial_{i} + \bar{\partial}_i$ are
tangent to $\Ga$.
\end{theo}

The proof starts from the following representation of the Schwartz
kernel of the Toeplitz operator $T_k$:
\begin{gather} \label{eq:noyau_Toeplitz}
 T_k (x,y) = \Bigl( \frac{k}{2 \pi} \Bigr)^n E^{k} (x,y) \tilde{f}
(x,y,k)  + O (k^{-\infty}) \end{gather}
where, if we consider $M^2$ as a complex manifold with complex
structure $(j , -j)$,
\begin{itemize}
\item 
$E$ is a section of $L \boxtimes \bar{L}
\rightarrow M^2$ satisfying  
$$ E (x,x) = u \otimes \bar{u}, \quad \forall u \in L_x \text{ such
  that }
  \| u \| = 1,
$$ $ \bar{\partial} E \equiv 0 $
modulo a section vanishing to every order along the diagonal $\De$ of
$M^2$ and  $\| E(x,y)
\| <1 $ if $x \neq y$.
\item
  $\tilde{f}(.,k)$ is a sequence of sections of $ \delta \boxtimes
  \bar{\delta} \rightarrow M^2$ with an asymptotic expansion of the form 
$$ \tilde{f}(.,k) = \tilde{f}_0 + k^{-1} \tilde{f}_1 + ...$$
whose coefficient satisfy $\bar{\partial} f_l \equiv 0  $
modulo a section vanishing to every order along $\De$.  Furthermore, $$\tilde{f}_0 (x,x) = f_0 (x),$$
where $f_0$ is the principal symbol of $T_k$. 
\end{itemize}
In other words, $T_k (.,.)$ is a
Lagrangian section associated to the diagonal $\De$ of $M^2$. This
result was proved in \cite{oim1}, without the additional bundle
$\delta$. The generalization is straightforward. 

Since the norm of $E$ is $<1$ outside the diagonal and $\Psi_k$ is
admissible, one has for every $x$ in $V$
$$ (T_k \Psi _k ) (x)= \Bigl( \frac{k}{2 \pi} \Bigr)^{n+ \frac{n}{4}}  \int E^k
(x,y).F^k(y) \ \tilde{f}(x,y,k). \tilde{g} (y,k)\ \mu_M(y) +
O(k^{-\infty})$$ 
where we integrate on a neighborhood of $x$. Introduce a unitary
section $t$ of $L$ over $U$ such that $\iota^* t = t_\Gamma$ over
$U_{\Gamma}$ and let us write  
\begin{gather} \label{eq:def_phi}
 E(x,y).F(y) = e^{i \phi(x,y) } t(x) .\end{gather} 
Then the imaginary part of $\phi$ is non positive and vanishes only if
$(x,y)$ belongs to 
$$ C:= \{ (x,x) \in M^2; \; x \in \iota(\Ga) \} .$$
To compute the derivatives of $\phi$ along $C$, recall the following
lemma proved in \cite{oim2} (cf. proposition 2.2 p.1535). 
\begin{lemme} \label{lem:derivee_F}
If $\nabla^{L} F = \tfrac{1}{i} \al_F \otimes F$, then
  $\al_F$ vanishes along $\iota(\Ga)$ and for every vector field $X,Y$ 
$$ \Lie_X \langle \al_F , Y \rangle  = \om ( q X, Y)$$
at $x \in \iota(\Ga)$, where $q$ is the projection onto $T_x^{0,1}M$ with kernel
$T_x \iota(\Ga) \otimes \C$. 
\end{lemme} 
As a corollary, we have the following
\begin{lemme} \label{lem:derivee_E}
If $\nabla^{L \boxtimes \bar{L}} E =
\tfrac{1}{i} \al_E \otimes E$, then $\al_E$ vanishes along the
diagonal, and for every vector fields $X_1,Y_1,X_2,Y_2$ of
$M$,
$$ \Lie_{(X_1,Y_1)}. \langle \al_E, (X_2,Y_2) \rangle = \om (
X_1^{0,1} - Y_1^{0,1}, X_2) + \om (
X_1^{1,0} - Y_1^{1,0}, Y_2)$$
along the diagonal, where we denoted by $X^{1,0}$ and $X^{0,1}$ the
holomorphic and anti-holomorphic parts of $X$ respectively.
\end{lemme}
We deduce from both lemmas that $d_y \phi$ vanishes along $C$. Furthermore
the kernel of the tangent map to $d_y \phi $ at $(x,x) \in C$ is 
$$ \bigl( T_x ^{0,1}M \times (0) \bigr) \oplus \bigl( T_{(x,x)} C \otimes \C \bigr) $$
Finally, we have along $C$, 
$$ d^2_y \phi  (Y_1,Y_2) = \om ( Y_1^{1,0}, Y_2) - \om ( q Y_1, Y_2) $$ 
and $d^2_y \phi $ is non-degenerate. So we can apply the stationary
phase lemma (cf. \cite{Ho} section 7.7 or theorem \ref{theo:PS1}).
One gets
$$ (T_k \Psi_k ) (x)   = \Bigl( \frac{k}{2 \pi} \Bigr)^{\frac{n}{4}} e^{ik
  \phi_r (x) } t^k(x) \ \tilde{h} (x,k) + O( k^{-\infty})$$
where 
\begin{gather} \label{eq:def_phir} 
\phi_r (x) \equiv \phi (x,y) \end{gather} 
modulo a linear combination with
  $\Ci$ coefficients of the
  $\partial_{y^i} \phi (x,y)$.
\begin{lemme} $e^{ik
  \phi_r (x) } t(x)$ satisfies the same assumption as the section
  $F$. 
\end{lemme}
\begin{proof} 
Since $\phi$ and $d_y \phi$ vanishes along $C$, $\phi_r$ vanishes
along $\iota(\Ga)$
and consequently $$ e^{i
  \phi_r (x) } t(x) = t_{\Ga} (x)$$ for every $x$ in
$\iota(\Ga)$. Introduce complex coordinates $x^1,...,x^n$ and write
$$ \nabla t = \tfrac{1}{i}  t \otimes \sum a_j (x) dx^j + \bar{a}_j
(x) d\bar{x}^j .$$
Derivating \eqref{eq:def_phi}, it follows from $\bar{\partial} E
\equiv 0$ that  
$$ \partial_{\bar{x}^i} \phi (x,y) \equiv \bar{a}_j (x) \mod
\ideal_{\De} (\infty),$$ i.e. modulo a function vanishing to infinite
order along the diagonal. Derivating again, one has
$$\partial_{\bar{x}^i} \partial_{y^j} \phi (x,y) \equiv 0 \mod
\ideal_{\De} (\infty) $$
Then we deduce from the two previous equations and
\eqref{eq:def_phir} that for every multi-index $\al$,
$$ \partial_{\bar{x}^1}^{\al(1)}...\partial_{\bar{x}^n}^{\al(n)} \bigl(
\partial_{\bar{x}^i} \phi_r  - \bar{a}_i \bigr) (x) = 0 $$
along $\iota(\Ga)$. And consequently 
$$ \partial_{\bar{x}^i} \phi_r  \equiv \bar{a}_i \mod \ideal_{\iota(\Ga)}
(\infty)$$
which proves the result. 
\end{proof}   

\begin{lemme} 
The sequence $\tilde{h}(.,k)$ admits an asymptotic
  expansion $\tilde{h}_0 + k^{-1} \tilde{h}_1 +...$  whose
  coefficients satisfy $\bar{\partial} \tilde{h}_\ell \equiv 0$ modulo a
  section vanishing to every order along $\Ga$. 
\end{lemme} 

\begin{proof} First one deduces from lemma \ref{lem:derivee_F} that
  the imaginary part of $\phi_r$ and its first derivatives vanishes
  along $\iota(\Ga)$. Furthermore the Hessian of $\Im \phi_r$ along $\iota(
  \Ga)$ is non-degenerate in the
  transverse direction to $\iota(\Ga)$. As a consequence, if $\tilde{e}
  (x,k)$ is a sequence with an asymptotic expansion $\tilde{e}_0 +
  k^{-1} \tilde{e}_1 +...$ such that 
$$ e^{ik \phi_r (x)}   \tilde{e} (x,k) = O ( k^{-\infty}) $$
 then the coefficients $\tilde{h}_\ell$ vanish to every order along
 $\iota(\Ga)$. This was proved in \cite{oim1} (cf. lemma 1, p.6). We apply this to the
 sequence $\bar{\partial} T_k \Psi_k$ which vanishes, since $ \Pi_k T_k
 = T_k$ implies that $T_k \Psi_k$ belongs to
$\Hilbert_k$.  
\end{proof} 

The two previous lemmas imply that $T_k \Psi_k$ is a Lagrangian section. 
Then applying theorems \ref{theo:PS1} and \ref{theo:PS2}, we obtain
the symbol of $T_k\Psi_k$ by computations of linear algebra, which are easily done using
the tangent vectors $\partial_i$ introduced in the statement of
theorem \ref{theo:main}.

\subsection{Proofs of the theorems of part \ref{sec:resultat}}
    
A first corollary of theorem \ref{theo:main} is the existence of a
Lagrangian section with an arbitrary symbol: applying theorem \ref{theo:main} with the Toeplitz operator $\Pi_k$,
we construct a Lagrangian section with a prescribed principal symbol,
then theorem \ref{existence} follows from Borel resummation. Theorem
\ref{Toep_Lag} is a particular case of theorem \ref{theo:main}.

To prove theorem \ref{resultat_clef}, we can assume that $f_0$
vanishes along $\Ga$. Since we compute the symbol of
$T_k \Psi_k$  modulo $O(\hb^2)$, we can replace $T_k$ with every
Toeplitz operators of symbol $f_0 + \hb f_1 + O(\hb^2)$. So by theorem
\ref{normalised_symbol_quantisation}, we can choose
$$ T_k = \Pi_k \Bigr( f_0 + k^{-1} f_1 + \frac{1}{ik} ( \nabla^{L^k}_X
\otimes \identite + \identite \otimes D_X ) \Bigl)  $$
where $X$ is the Hamiltonian vector field of $f_0$. So $T_k \Psi_k$ is
equal to 
\begin{gather} \label{La}
 \Pi_k \Bigl[ \bigl( \tfrac{k}{2 \pi} \bigr)^{\frac{n}{4}} \bigl( [f_0 + \tfrac{1}{ik} \nabla_{X}^{L^k} ] F^k
\bigr) \tilde{g}(.,k) \Bigr] + k^{-1} \Pi_k \Bigl[  \bigl( \tfrac{k}{2
  \pi} \bigr)^{\frac{n}{4}} F^k \bigl( [f_1 +
\tfrac{1}{i} D_X ] \tilde{g}(.,k) \bigr)  \Bigr] \end{gather}  
By theorem \ref{theo:main}, each term of the sum is a Lagrangian
section. Furthermore by the first part of this theorem, the symbol of the second one is 
$$ \hb\ \bigl( (\iota^* f_1).g_0 + \iota^*(\tfrac{1}{i} D_X
\tilde{g}_0) \bigr) +O(\hb^2).$$
Since $X$ is tangent to $\iota(\Ga)$, $\iota ^* (D_X
\tilde{g}_0)$ only depends on the restriction of $\tilde{g}_0$
to $\iota(\Ga)$. So we can define the operator $\iota^* D_X$ acting on $\Ci
(\Ga, \delta_\Ga)$ which sends $g_0$ to $ \iota ^* (D_X
\tilde{g}_0)$. And the previous symbol is 
\begin{gather} \label{2}
 \hb\ \bigl( (\iota^* f_1) + \tfrac{1}{i} (\iota^* D_X) \bigr).g_0
 +O(\hb^2). 
\end{gather}

To compute the symbol of the first term of \eqref{La}, let us write
$$  [f_0 + \tfrac{1}{ik} \nabla_{X}^{L^k} ] F^k
 = F^k a $$
with $a$ defined on a neighborhood of $\iota(\Ga)$.  
\begin{lemme} \label{li2}
The function $a$ and its first derivatives vanish along $\iota(\Ga)$. If $Z$
and $W$ are holomorphic vector fields of $M$, then
$$ \Lie_{\bar{W}} \Lie_{\bar{Z}} a =  \om( \bar{W} , \Lie_X \bar{Z}  ) $$
on $\iota( \Ga)$.
\end{lemme}

\begin{proof} 
Denote by $\alpha _F$ be the one-form such that $ \nabla^{L} F = 
\frac{1}{i}\al_F \otimes F$. Then
$$ a = f_0 - \langle \al_{F} , X \rangle .$$
This vanishes along $\iota ( \Ga)$ because $\al_F$ vanishes along $\iota(\Ga)$
(cf. lemma \ref{lem:derivee_F}).

Since $X$ is the Hamiltonian vector field of $f_0$, one has $\Lie_Y
f_0 + \om (X,Y) =0$. Since the curvature of $L$ is $\frac{1}{i}\om$,
one has $d \al_F = \om$ and consequently 
$$\Lie_Y \langle \al_F, X
\rangle = \Lie_X \langle \al_F, Y
\rangle + \om(Y,X) + \langle \al_F, [Y,X] \rangle.$$ It follows that
$$ \Lie_Y a = - \Lie_X \langle \al_F, Y \rangle - \langle \al_F, [Y,X]
\rangle .$$
This vanishes along $\iota(\Ga)$, because $\al_F$ vanishes along $\iota(\Ga)$
and $X$ is tangent to $\iota(\Ga)$. 

Since $\bar{\partial} F \equiv 0$ modulo a flat section and $Z$ is holomorphic, $  \langle
\al_F, \bar{Z} \rangle$ vanishes to every order along $\iota(\Ga)$. So choosing
$Y =\bar{Z}$ in the previous equation, we obtain
$$ \Lie_{\bar{W}}.\Lie_{\bar{Z}} a = - \Lie_{\bar{W}} \langle \al_F,
[\bar{Z}, X] \rangle \qquad \text{ along } \iota(\Ga) .$$
Using again that $\om = d\al_F$ and $\al_F $ vanishes along $\iota (\Ga)$, it
follows that
$$  \Lie_{\bar{W}}.\Lie_{\bar{Z}} a = - \om ( \bar{W}, [\bar{Z}, X ] )
- \Lie_{[\bar{Z},X ]} \langle \al_F, \bar{W} \rangle $$
along $\iota ( \Ga)$. The second term of the right side vanishes along
$\iota ( \Ga)$ because $\bar{W}$ is an anti-holomorphic vector
field. This gives the result.
\end{proof} 

Since $\iota^* D_X$ and $\Lie_X^{\delta_{\Ga}}$ are first order
differential operators which have the same symbol,
$$ \iota^* D_X - \Lie_X^{\delta_{\Ga}} = b $$
where $b$ is a function of $\Ga$. 

\begin{lemme} \label{li3} 
The symbol of $ \Pi_k \bigl[ F^k a \tilde{g}(.,k) \bigr]$ is 
$  i  \hb  b. g_0 + O(\hb^2)$
\end{lemme}

So the symbol of $T_k \Psi_k$ is the sum of \eqref{2} and $ i 
\hb b g_0 + O(\hb^2)$, which is equal to 
$$ \hb(\iota^* f_1 +  \tfrac{1}{i} \Lie_X^{\delta_{\Gamma}}) . g_0 + O(\hb^2). $$
Theorem \ref{resultat_clef} follows. 

\begin{proof} 
Let us start with a local computation of the function $b$.  
Let $u$ be a non-vanishing section of $\delta_{\Ga} \rightarrow
\Ga$. Then one has 
$$  b = \frac{ \bigl( \iota^* D_X - \Lie_X^{\delta_{\Ga}} \bigr) .u
}{u}$$
Let $\beta$ be a non-vanishing $(n,0)$-form of $M$ such that $\varphi
(u^{\otimes 2}(x)) = \beta (\iota(x))$ if $x$ be\-longs to $\Ga$. Then
$\iota^* D_X $ is defined in such a way that 
$$ \frac{ ( \iota^* D_X).u}{u} = \frac{1}{2} \iota^* \Bigl(\frac{ p
  \Lie_{X} \be}{\be}\Bigr) = \frac{1}{2} \frac{ \iota^* (p
  \Lie_{X} \be)}{\iota^* \be} .$$
where $p$ is the projection from $\La^{n} M \otimes \C$ onto
  $\La^{n,0} M$ with kernel 
$ \La^{0,n} M \oplus ... \oplus \La^{n-1,1} M.$ 
On the other hand, since 
$\varphi_{\Ga}(u^{\otimes 2} (x)) = \iota^* \varphi
  (u^{\otimes 2} (x) ) = \iota^* \be (x) $, 
one has
$$ 
 \frac{ \Lie_X^{\delta_{\Ga}}.u
}{u} = \frac{1}{2} \frac{ \Lie_X \iota^* \be}{\iota^* \be} = \frac{1}{2} \frac{ \iota^* \Lie_X \be}{\iota^* \be} .$$
Consequently
$$ b = \frac{1}{2} \frac{ \iota^* \bigl(p \Lie_X \be-  \Lie_X
  \be\bigr) }{\iota ^* \be}. $$

Now let us choose a frame $(\partial_1,...,\partial_n)$ of holomorphic
  vector fields of $M$ such that the vectors $\partial_i + \bar{\partial_i}$ are tangent
  to $\iota(\Ga)$. Denote by $(\te^{1},..., \te^n)$ the dual frame and set
$$ \be = \te^1 \wedge... \wedge \te^n .$$ 
Then using that $\Lie_X \te^i \equiv - \sum \langle \te^i, \Lie_X
\bar{\partial}_j \rangle \bar{\te}^j $ modulo a linear combination of the
$\te^i$, we obtain that $\iota^* (p \Lie_X \be-  \Lie_X
  \be)$ is equal to 
\begin{gather*} 
\sum \langle \te^1, \Lie_X \bar{\partial}_j \rangle
  \bar{\te}^j \wedge \te^2 \wedge ...\wedge \te^n +  \langle \te^2,
  \Lie_X \bar{\partial}_j \rangle \te^1 \wedge \bar{\te}^j \wedge
  \te^3 \wedge ...\wedge \te^n + ...\\
+  \langle \te^n,
  \Lie_X \bar{\partial}_j \rangle \te^1 \wedge ...\wedge \te^{n-1}
  \wedge \bar{\te}^j \end{gather*} 
It follows then from $\iota^* \te^i  = \iota^* \bar{\te}^i$ that 
$$ b = \frac{1}{2} \sum \langle \te^i, \Lie_X
\bar{\partial}_i \rangle .$$

To end the proof, assume furthermore that $\frac{1}{i} \om ( \partial_i
,\bar{\partial}_j) = \de_{ij}$. Then it follows from the second part of theorem \ref{theo:main}
that the symbol of $ \Pi_k \bigl[ F^k a \tilde{g}(.,k) \bigr]$ is 
$$ - \hb g_0  \frac{1}{2} \ \iota^* \sum \Lie_{\bar{\partial}_i}
\Lie_{\bar{\partial}_i} a $$
which by lemma \ref{li2} is equal to  
$$ - \hb g_0 \frac{1}{2}  \ \iota^* \sum \om ( \bar{\partial}_i, \Lie_X
\bar{\partial}_i)$$
Using again that $\frac{1}{i} \om ( \partial_i
,\bar{\partial}_j) = \de_{ij}$, we obtain
 $$   \hb g_0 \frac{i}{2} \ \iota^*  \sum \langle \te^i, \Lie_X
\bar{\partial}_i \rangle .$$
The final result follows. 
\end{proof} 

Finally, let us prove theorem \ref{norme}.
One has
$$ \int_M \xi \  \bigl\| \Psi_k \bigr\|_{L^k \otimes \delta} ^2 \; \mu_M =
\Bigl( \frac{k}{2\pi} \Bigr)^{\frac{n}{2}} \int_M  e^{-k c} \ \xi \
\bigl\| \tilde{g}(.,k) \bigr\|_{\delta}^2 \; 
\mu_M +O(k^{-\infty}) $$
where $c(x) = -2 \ln \| F(x) \| _{L}$. The following is a
consequence of lemma
\ref{lem:derivee_F}.
\begin{lemme} \label{lem:norme_F}
The function $c$ and its first derivatives vanish along
$\iota(\Ga)$. Furthermore, its Hessian at $x \in \iota(\Ga)$ is definite positive on
$J T_x \iota(\Ga)$ and is given by 
$$  X.Y.c  = 2 \om (X, JY), \quad X,Y \in J T_x \iota(\Ga).$$
\end{lemme} 
So integrating along transversal directions to $\Ga$, it follows from the
stationary phase lemma that 
$$  \int_M \; \xi \ \bigl\| \Psi_k \bigr\|_{L^k  \otimes \delta}
^2 \;  \mu_M = \int_{\Gamma} (\iota^* \xi)\  d +O(k^{-1})$$
where $d$ is the density of $\Ga$ such that 
$$ d\eval{x} ( X) =   \bigl\| g_0(x) \bigr\|^2_{ \delta}
\  \mu_{M} \eval{x} (X \wedge Y) \ 2^{-\frac{n}{2}} \Bigl( \det \bigl[ \om (Y_i,J Y_j)
\bigr]\Bigr)^{-\frac{1}{2}}.$$
Here $(X_i)$ and $(Y_i)$ are bases of $T_x \Ga$ and $J T_x \Ga$
respectively, and $X = X_1 \wedge...\wedge X_n$, $Y = Y_1
\wedge ... \wedge Y_n$. To deduce
theorem \ref{norme}, we have to check that  
\begin{gather*}
d = m (g_0)  
\end{gather*}
This is easily done, by introducing the same vector fields
$\partial_i$ and forms $\te^i$ like in the
proof of lemma \ref{li3}, setting $X_i = \frac{1}{\sqrt{2}}(
\partial_i + \bar{\partial}_i)$, $Y_j = J X_j$ and choosing $g_0$ such that $ \varphi (g_0^{\otimes
  2} ) = \theta_1 \wedge ...\wedge \theta^n$. 

\scratch{Quelques d{\'e}tails: Pour d{\'e}duire le lemme \ref{lem:norme_F}
  du lemme \ref{lem:derivee_F}, on utilise que si $X \in J
  T\iota(\Ga)$, alors $q(JX) = 0$ et $q(X+ i JX) = X+iJX$, donc $q(X)
  = X + iJX$. Apr{\`e}s, on calcule
$$  X.Y .c = - \frac{X .Y. \| F \|^2}{\| F \|^2} = - \frac{1}{i} X. \langle \al_F -
\bar{\al}_F, Y \rangle .$$  
Ce qui donne avec le lemme \ref{lem:derivee_F}
$$ -\frac{1}{i} \om (X + iJX - (X - iJX), Y) 
= -2 \om ( JX,Y) = 
2 \om (X, JY).$$
Voil{\`a}. Ensuite, pour la v{\'e}rification de $d = m(g_0)$, avec les
choix que j'indique $\om (Y_i, J Y_j) =\de_{ij}$, $ \mu_M ( X \wedge
Y) = 1$, et $\te^1 \wedge ...\wedge \te^n$ {\'e}tant de norme 1, il
vient $d(X) = 2^{-\frac{n}{2}}$. D'autre part $\langle \iota^*
\te^i, X_j \rangle = \frac{1}{\sqrt{2}} \de_{ij} $, donc $\langle \iota^*
\varphi (g_0^2) , X \rangle =  2^{-\frac{n}{2}}$, autrement dit $m(
g_0) (X) = 2^{- \frac{n}{2}}$.  
}

\section{Appendix} 

Let $W$ be an open set of $\R^{n} \times \R^{k} \ni (x,y)$. Denote by
$p$ the projection from $W$ onto $\R^n$. Let $\ph (x,y)$ be a $\Ci$
function on $W$ whose imaginary part is $ \geqslant 0$. Let $a(x,y)$
be a $\Ci$ function with compact support in $W$. Stationary phase
lemma gives the asymptotic expansion  of  
$$ I(a, \ph)(x, \tau ) = \int_{\R^k} e^{i \tau \ph (x,y) } a(x,y) |dy| $$
when $\tau \rightarrow \infty$.
First, if the support of $a$ doesn't meet the critical locus 
$$ C := \{ (x,y)\in W  ; \; d_y \ph (x,y) = 0 \text{ and } \Im \ph (x,y) = 0 \}, $$ 
then  $I(a, \ph)$ is $O(\tau^{-\infty})$. Introduce the functions
$$\ph_{i,j}(x,y) = \partial_{y^i} \partial_{y^j} \ph
(x,y), \quad i,j =1,...,k.$$ 
The following theorem is proved in \cite{Ho} section 7.7. 

\begin{theo} \label{theo:PS1}
Assume that at $(x_0,y_0) \in C$, the matrix $
\bigl(\ph_{i,j}(x_0,y_0)\bigr)$ is invertible. Then there exists a
neighborhood $U$ of $(x_0, y_0)$ such that if the support of $a$ is a
subset of $ U$, one has 
$$ I(\ph, a) (x, \tau) =( \tfrac{2 \pi }{ \tau} )^{\frac{k}{2}} \ d (x) \
e^{i \tau \ph_r (x) }\ b(x,\tau) + O(\tau^{-\infty}) \qquad \text{ over
} p(U)
$$ 
where  $d$, $\ph_r$ et $b(.,\tau)$ are $\Ci$ functions such that  
\begin{itemize} 
\item $d$ only depends on $\ph$. In particular, 
 $$ d(x) = \operatorname{det}^{-\frac{1}{2}} [ \tfrac{1}{i}
 \varphi_{j,k} (x,y)]_{j,k}, \quad  \text {if } (x,y) \in C \cap U.$$
\item $\ph_r$ is such that $ \ph (x,y) \equiv \ph_r (x) $ on $U$
  modulo a linear combination with $\Ci$ coefficient of the functions  $\partial_{y^j} \ph$. 
\item  $b(.,\tau)$ has an asymptotic expansion for the $\Ci$ topology of the form 
  $$ b_0 (x) + \tau b_1(x) + \tau^2 b_2 (x) +...$$ 
Furthermore,   $b_0(x) = a(x,y)$
if  $(x,y) \in C \cap U$. 
\end{itemize} 
\end{theo} 
In \cite{Ho} the various terms of the asymptotic expansion are
completely determined and not only their restriction at $p(C)$. But in the applications this leads to
complicated computations that we prefer to avoid. Let us introduce an
additional assumption.

Denote by $E_{(x,y)}$ the complexification of the tangent space to the
fiber of $p$ at  $(x,y)$. 
At $(x,y) \in C$ the tangent map to the section  $d_y
\ph$ of $E^*$ 
$$ T_{(x,y)} d_{y} \ph : T_{(x,y)} W \otimes \C \rightarrow E^*_{(x,y)} $$ 
 is well-defined. Assume that  $(\varphi_{i,j})$ is invertible along
 $C$, that is the kernel   $F_{(x,y)}$ of $ T_{(x,y)} d_{y} \ph$
satisfies  
\begin{gather} \label{H1}
 \forall \ (x,y) \in C, \quad  F_{(x,y)} \oplus E_{(x,y)}  =
  T_{(x,y)} W \otimes \C. \end{gather}
Assume furthermore that 
\begin{gather} \label{H2}
\text{$C$ is a submanifold of $W$ and } 
T C \otimes \C = F \cap
  \bar{F} . \end{gather} 
Finally, these two assumptions imply that the restriction $p : C \rightarrow
  \R^n$ is an immersion. We assume it is an embedding. 
   
Observe that when the phase takes real values, the assumption
\eqref{H2} is a consequence of 
\eqref{H1}. We are interested in the opposite case, typically when the
Hessian of the imaginary part of the phase is non-degenerate  in the
transverse directions to $C$, for instance with $$\ph (x,y) = xy + \tfrac{i}{2}(x^2 +
y^2).$$ 
We can also consider intermediary cases, for example  $\ph (x,y) = xy + \frac{i}{2} y^2$. 

Under the previous assumptions, when the amplitude $a$ vanishes to
order $m$ along $C$, i.e. when the partial derivatives of $a$ of order 
$\leqslant m-1$ vanish along $C$, it follows from the result of
\cite{Ho} that the functions $b_i$ vanish to order $m-2i$ along
$C$. Furthermore one can easily compute 
$b_i$ modulo a function vanishing to order   $m-2i+1$ along 
$C$. 

To state the result, consider a free family $\partial_1,...,
\partial_l$ of complex tangent vectors to $W$ at  $(x,y) \in
   C$ such that 
$$ \operatorname{Vect}_{\C} (\partial_1,..., \partial_l) \oplus  (TC
\otimes \C) = F_{(x,y)}.$$
If $a$ vanishes to order  $m$ along  $C$, we define the polynomial 
$$ [a] (Z,Y) = \sum_{|\al| + |\be|=m} \frac{1}{\al ! \be !} \Bigl( 
\partial_1^{\al(1)}...\partial_l^{\al(l)}\partial_{y^1}^{\be(1)}...\partial_{y^k}^{\be(l)}
a(x,y) \Bigr) \ Z^{\al} Y^{\be}$$
at $(x,y) \in C$. 
Similarly, if $b(x)$ vanishes to order  $l$ along  $p(C)$, we set 
$$ [b] (Z) = \sum_{|\al|=l} \frac{1}{\al ! \be !} \Bigl( 
(p_*\partial_1)^{\al(1)}...(p_* \partial_l)^{\al(l)} b(x) \Bigr) \ Z^{\al}.$$
at $x \in p(C)$. 
\begin{theo} \label{theo:PS2}
Under the assumptions \eqref{H1} and \eqref{H2}, if $a$ vanishes to
  order  $m$ along  $C$, then for every  $i \leqslant \frac{m}{2}$ the
  function  $b_i$ vanishes to order  $m-2i$ along  $p(C)$. Furthermore  
$$ [b_i](Z) = \frac{1}{i!} \De^i A_{2i}(Z,Y) $$
at $(x,y) \in C$, where 
\begin{itemize} 
\item $[a](Z,Y) = \sum_{l=0}^{m} A_l(Z,Y)$ and $A_l$ is homogeneous of 
 degree $l$ in $Y$ and of degree  $m-l$ in $Z$. 
\item  $\De = \tfrac{i}{2} \textstyle{\sum}_{j,k}
  \ph^{j,k}(x,y) \partial_{Y^j} \partial_{Y^k}$, with  $(\ph^{j,k}(x,y))$
  the inverse of $ (\ph_{j,k}(x,y))$.
\end{itemize} 
\end{theo} 

More intrinsically, denote by $\ideal ^m(C) \subset \Ci( W)$ the ideal
of functions vanishing to order  $l$ along  $C$. Then $\ideal ^m (C)/
\ideal ^{m-1} (C)$ is isomorphic to the space of sections of the $m$-th symmetric power of the
complex conormal bundle $\no^* (C)$. By \eqref{H1}, we have an
isomorphism of vector bundle over $C$, 
$$  \no^* (C) \simeq \no ^* ( p (C))  \oplus E^* ,$$
which associates  $\ga$ and $\al \oplus \be$ if 
$$ \langle \ga, U +
V \rangle = \langle \al, p_* U \rangle + \langle \beta, V \rangle,
\quad \text{$U \in F$ and $V \in  E$} $$
Consequently,  
\begin{gather*} 
 \Sym_m (  \no^* (C) ) = \bigoplus_{l=0}^{m} \Sym_{m-l} \bigl( \no^* (p(C) )  \bigr) \otimes \Sym_{l} ( E^* ) \end{gather*} 
$ \Delta = \tfrac{i}{2} \textstyle{\sum}_{j,k} \ph^{j,k}
\partial_{y^j} \partial_{y^k}$ defines a section of $\Sym_{2}(E)$, so
$\Delta ^i$ acts as an operator 
$$\Sym_{2i} (
E^* ) \rightarrow \C  .$$ In theorem  \ref{theo:PS2}, we consider
$[a]$ as a section of $\Sym_m (  \no^*
(C) )$ and  $[b_i]$ as a section of  $\Sym_{m-2i} ( \no^* (p(C)))$ . Then we have  
$$ [b_i] = \frac{1}{i!} (\identite \otimes \De^i) A_{2i} $$
where 
$$ [a] = \sum_{l=0}^{m} A_l, \quad A_l \in \Ci (C,\Sym_{m-l} \bigl( \no^*
(p(C) )  \bigr) \otimes \Sym_{l} ( E^* )) .$$

\bibliography{biblio}

\begin{thebibliography}{10}

\bibitem{Be}
F.~A. Berezin.
\newblock General concept of quantization.
\newblock {\em Comm. Math. Phys.}, 40:153--174, 1975.

\bibitem{BeCaGu}
M{\'e}lanie Bertelson, Michel Cahen, and Simone Gutt.
\newblock Equivalence of star products.
\newblock {\em Classical Quantum Gravity}, 14(1A):A93--A107, 1997.
\newblock Geometry and physics.

\bibitem{BoMeSc}
M.~Bordemann, E.~Meinrenken, and M.~Schlichenmaier.
\newblock Toeplitz quantization of {K}\"ahler manifolds and ${\rm gl}({N})$,
  ${N}\to\infty$ limits.
\newblock {\em Comm. Math. Phys.}, 165(2):281--296, 1994.

\bibitem{BP2}
D.~Borthwick, T.~Paul, and A.~Uribe.
\newblock Semiclassical spectral estimates for {T}oeplitz operators.
\newblock {\em Ann. Inst. Fourier (Grenoble)}, 48(4):1189--1229, 1998.

\bibitem{BoGu}
L.~Boutet~de Monvel and V.~Guillemin.
\newblock {\em The spectral theory of {T}oeplitz operators}, volume~99 of {\em
  Annals of Mathematics Studies}.
\newblock Princeton University Press, Princeton, NJ, 1981.

\bibitem{oim1}
L.~Charles.
\newblock Berezin-{T}oeplitz operators, a semi-classical approach.
\newblock {\em Comm. Math. Phys.}, 239(1-2):1--28, 2003.

\bibitem{oim2}
L.~Charles.
\newblock Quasimodes and {B}ohr-{S}ommerfeld conditions for the {T}oeplitz
  operators.
\newblock {\em Comm. Partial Differential Equations}, 28(9-10):1527--1566,
  2003.

\bibitem{new}
L.~Charles.
\newblock Semi-classical properties of geometric quantization with metaplectic
  correction.
\newblock http://www.institut.math.jussieu.fr/\~{}charles/Articles/Half2.pdf,
  2006.

\bibitem{Co}
Y.~Colin~de Verdi\`eres.
\newblock Bohr-sommerfeld rules to all orders.
\newblock to appear in {\em Ann. Henri Poincar\'e},
  http://www-fourier.ujf-grenoble.fr/\~{}ycolver/ebk.ps, 2003.

\bibitem{Du3}
J.~J. Duistermaat.
\newblock Oscillatory integrals, {L}agrange immersions and unfolding of
  singularities.
\newblock {\em Comm. Pure Appl. Math.}, 27:207--281, 1974.

\bibitem{Fe}
Boris Fedosov.
\newblock {\em Deformation quantization and index theory}, volume~9 of {\em
  Mathematical Topics}.
\newblock Akademie Verlag, Berlin, 1996.

\bibitem{Gu}
Victor Guillemin.
\newblock Star products on compact pre-quantizable symplectic manifolds.
\newblock {\em Lett. Math. Phys.}, 35(1):85--89, 1995.

\bibitem{Ho}
Lars H{\"o}rmander.
\newblock {\em The analysis of linear partial differential operators. {I}},
  volume 256 of {\em Grundlehren der Mathematischen Wissenschaften [Fundamental
  Principles of Mathematical Sciences]}.
\newblock Springer-Verlag, Berlin, second edition, 1990.
\newblock Distribution theory and Fourier analysis.

\bibitem{ka}
A.~V. Karabegov.
\newblock Cohomological classification of deformation quantizations with
  separation of variables.
\newblock {\em Lett. Math. Phys.}, 43(4):347--357, 1998.

\bibitem{KaSc}
Alexander~V. Karabegov and Martin Schlichenmaier.
\newblock Identification of {B}erezin-{T}oeplitz deformation quantization.
\newblock {\em J. Reine Angew. Math.}, 540:49--76, 2001.

\bibitem{Ko}
Bertram Kostant.
\newblock Quantization and unitary representations. {I}. {P}requantization.
\newblock In {\em Lectures in modern analysis and applications, III}, pages
  87--208. Lecture Notes in Math., Vol. 170. Springer, Berlin, 1970.

\bibitem{McSa}
Dusa McDuff and Dietmar Salamon.
\newblock {\em Introduction to symplectic topology}.
\newblock Oxford Mathematical Monographs. The Clarendon Press Oxford University
  Press, New York, second edition, 1998.

\bibitem{NeTs}
Ryszard Nest and Boris Tsygan.
\newblock Algebraic index theorem.
\newblock {\em Comm. Math. Phys.}, 172(2):223--262, 1995.

\bibitem{So}
J.-M. Souriau.
\newblock {\em Structure des syst\`emes dynamiques}.
\newblock Ma\^itrises de math\'ematiques. Dunod, Paris, 1970.

\bibitem{Tu}
G.~M. Tuynman.
\newblock Quantization: towards a comparison between methods.
\newblock {\em J. Math. Phys.}, 28(12):2829--2840, 1987.

\bibitem{We}
Alan Weinstein.
\newblock The {M}aslov gerbe.
\newblock {\em Lett. Math. Phys.}, 69:3--9, 2004.

\bibitem{Ze}
Steven Zelditch.
\newblock Index and dynamics of quantized contact transformations.
\newblock {\em Ann. Inst. Fourier (Grenoble)}, 47(1):305--363, 1997.

\end{thebibliography}

\end{document}